\documentclass{amsproc}
\usepackage{breqn,float,amssymb,pdflscape,rotating,hyperref}
\makeatletter
\@namedef{subjclassname@2020}{%
  \textup{2020} Mathematics Subject Classification}
\makeatother
\newtheorem{theorem}{Theorem}[section]

\theoremstyle{definition}

\newtheorem{example}[theorem]{Example}

\theoremstyle{remark}

\numberwithin{equation}{section}
\allowdisplaybreaks
\begin{document}

\title{An extended Cauchy integral}


\author{Robert Reynolds}
\address[Robert Reynolds]{Department of Mathematics and Statistics, York University, Toronto, ON, Canada, M3J1P3}
\email[Corresponding author]{milver73@gmail.com}
\thanks{}




\date{}

\dedicatory{}

\begin{abstract}
A new integral representation is derived using a definite integral given by Cauchy and used to evaluate a number of integrals containing the finite series of special functions.
\end{abstract}

\maketitle
\section{Introduction}
In the book \cite{bdh} by Bierens de Haan on definite integrals the following interesting integral formula is given:
\begin{multline}\label{eq1}
\int_0^{\infty } \frac{x^{p-1}}{\left(1+x^a\right) \left(1+x^b\right)} \, dx=\frac{\pi  }{2 a \sin (p \pi )}\sum _{n=0}^{a-1}
   \frac{\cos \left(\frac{(2 n-a+1) p \pi }{a}\right)+\cos \left(\frac{(2 n-a+1) (p-b) \pi }{a}\right)}{1+\cos
   \left(\frac{(2 n+1) (b \pi )}{a}\right)}\\
+\frac{\pi  }{2 b \sin (p \pi )}\sum _{n=0}^{b-1} \frac{\cos \left(\frac{(2
   n-b+1) p \pi }{b}\right)+\cos \left(\frac{(2 n-b+1) (p-a) \pi }{b}\right)}{1+\cos \left(\frac{(2 n+1) (a \pi
   )}{b}\right)}
\end{multline}
where $Re(p)>0,a\neq b$ and $a,b$ are even numbers. The reference given for this formula is equation (59) in Verhandelingen der Koninklijke Akademie van Wetenschappen, Deel IV, 1858., Tables d'int\'{e}grales d\'{e}fines, par D. Bierens de Haan, 
[\href{https://archive.org/details/verhandelingende04koni/page/n9/mode/2up}{BDH}]. Following this reference leads to the work by Cauchy on definite integrals published in 1827, M\'{e}m. Paris. 1823. 603. Sav. Etr. I. 1827. p.3, [\href{http://sites.mathdoc.fr/cgi-bin/oeitem?id=OE_CAUCHY_1_1_319_0}{Cauchy}]. These work were used in the publishing of the work by Gradshteyn and Ryzhik \cite{grad}, et al. Prudnikov et al, \cite{prud1} and Gr\"{o}bner \cite{grobner}, to name a few. The format of these  very useful tables of definite integrals are a detailed list of integrals involving special functions, finite and infinite series and algebraic functions without proofs. In this current work we produce a proof and evaluations of the extended form of the Cauchy integral. The proof involves the contour integral method in \cite{reyn4} and the evaluations involve the finite series of the Lerch-Hurwitz zeta function [DLMF,\href{https://dlmf.nist.gov/25.14}{14}].
\subsection{Contour integral representations}
We proceed by applying the contour integral method in \cite{reyn4} to equation (\ref{eq1}) to get
\begin{multline}\label{eq2}
\frac{1}{2\pi i}\int_{C}\frac{a^w w^{-k-1} x^{m+w-1}}{\left(x^{\alpha }+1\right) \left(x^{\beta }+1\right)}dw
=\frac{1}{2\pi i}\int_{C}\sum_{j=0}^{\beta-1}\frac{\pi  a^w w^{-k-1} \csc
   (\pi  (m+w)) \cos \left(\frac{\pi  (-\beta +2 j+1) (m+w)}{\beta }\right)}{2 \beta  \left(\cos \left(\frac{\pi 
   \alpha  (2 j+1)}{\beta }\right)+1\right)}dw\\
+\frac{1}{2\pi i}\int_{C}\sum_{j=0}^{\beta-1}\frac{\pi  a^w w^{-k-1} \csc (\pi  (m+w)) \cos \left(\frac{\pi  (-\beta +2
   j+1) (-\alpha +m+w)}{\beta }\right)}{2 \beta  \left(\cos \left(\frac{\pi  \alpha  (2 j+1)}{\beta
   }\right)+1\right)}dw\\
+\frac{1}{2\pi i}\int_{C}\sum_{j=0}^{\alpha-1}\frac{\pi  a^w w^{-k-1} \csc (\pi  (m+w)) \cos \left(\frac{\pi  (-\alpha +2 j+1) (m+w)}{\alpha
   }\right)}{2 \alpha  \left(\cos \left(\frac{\pi  \beta  (2 j+1)}{\alpha }\right)+1\right)}dw\\
+\frac{1}{2\pi i}\int_{C}\sum_{j=0}^{\alpha-1}\frac{\pi  a^w w^{-k-1}
   \csc (\pi  (m+w)) \cos \left(\frac{\pi  (-\alpha +2 j+1) (-\beta +m+w)}{\alpha }\right)}{2 \alpha  \left(\cos
   \left(\frac{\pi  \beta  (2 j+1)}{\alpha }\right)+1\right)}dw.
\end{multline}
Next using \cite{reyn4} and a few algebraic processes, we form the contour integral representation involving the product of cosine-cosecant functions which is used for deriving the right-hand side of equation (\ref{eq2}) given by;
\begin{multline}\label{eq3}
-\frac{i (2 i)^k b^k e^{i (b m-(m+t) x)} \left(\Phi \left(e^{2 i b m},-k,\frac{b-x-i \log (a)}{2 b}\right)+e^{2
   i (m+t) x} \Phi \left(e^{2 i b m},-k,\frac{b+x-i \log (a)}{2 b}\right)\right)}{\Gamma(k+1)}\\
   =\frac{1}{2\pi i}\int_{C}a^w w^{-1-k} \cos ((m+t+w) x)
   \csc (b (m+w))dw.
\end{multline}
Here we again use the method in \cite{reyn4} to derive the left-hand side contour integral representation for equation (\ref{eq1}) given by;
\begin{equation}\label{eq4}
\int_{0}^{\infty}\frac{x^{m-1} \log ^k(a x)}{\Gamma(k+1) \left(x^{\alpha }+1\right) \left(x^{\beta }+1\right)}dx=\int_{0}^{\infty}\int_{C}\frac{a^w w^{-k-1}
   x^{m+w-1}}{\left(x^{\alpha }+1\right) \left(x^{\beta }+1\right)}dwdx
\end{equation}
We are now in a good position to write down the contour integral representation for equation (\ref{eq1}) using equations (\ref{eq3}) and (\ref{eq4}). To derive the first contour integral representation on the right-hand side of equation (\ref{eq2}) we set use equation (\ref{eq3}) and set $b= \pi ,t= 0,x= \frac{\pi  (-\beta +2 j+1)}{\beta }$ multiply both sides by $\frac{\pi }{2 \beta  \left(\cos \left(\frac{\pi  \alpha  (2 j+1)}{\beta }\right)+1\right)}$ and simplify to get;
\begin{multline}\label{eq5}
-\sum_{j=0}^{\beta-1}\frac{2^{k-1} (i \pi )^{k+1} e^{i \left(\pi  m-\frac{\pi  m (-\beta +2 j+1)}{\beta }\right)} }{\beta  \Gamma(k+1) \left(\cos \left(\frac{\pi  \alpha  (2 j+1)}{\beta }\right)+1\right)}\left(\Phi \left(e^{2 i m \pi },-k,\frac{-\frac{\pi  (2 j-\beta +1)}{\beta }-i \log (a)+\pi }{2 \pi }\right)\right. \\ \left.
+e^{\frac{2 i \pi m (-\beta +2 j+1)}{\beta }} \Phi \left(e^{2 i m \pi },-k,\frac{\frac{\pi  (2 j-\beta +1)}{\beta }-i \log (a)+\pi }{2\pi }\right)\right)\\
=\frac{1}{2\pi i}\int_{C}\sum_{j=0}^{\beta-1}\frac{\pi  a^w
   w^{-k-1} \csc (\pi  (m+w)) \cos \left(\frac{\pi  (-\beta +2 j+1) (m+w)}{\beta }\right)}{2 \beta  \left(\cos
   \left(\frac{\pi  \alpha  (2 j+1)}{\beta }\right)+1\right)}dw.
\end{multline}
Next we look at deriving the third contour integral representation on the right-hand side of equation (\ref{eq2}), by using equation (\ref{eq3}) and setting $b= \pi ,t= -\alpha ,x= \frac{\pi  (-\beta +2 j+1)}{\beta }$ multiplying both sides by $\frac{\pi }{2 \beta  \left(\cos \left(\frac{\pi  \alpha  (2 j+1)}{\beta }\right)+1\right)}$ and simplifying to get;
\begin{multline}\label{eq6}
-\sum_{j=0}^{\beta-1}\frac{2^{k-1} (i \pi )^{k+1} e^{i \left(\pi  m-\frac{\pi  (-\beta +2 j+1) (m-\alpha )}{\beta }\right)}
   }{\beta  \Gamma(k+1) \left(\cos \left(\frac{\pi  \alpha  (2 j+1)}{\beta
   }\right)+1\right)}\left(e^{\frac{2 i \pi  (-\beta +2 j+1) (m-\alpha )}{\beta }} \Phi \left(e^{2 i m \pi },-k,\frac{\frac{\pi  (2
   j-\beta +1)}{\beta }-i \log (a)+\pi }{2 \pi }\right)\right. \\ \left.
+\Phi \left(e^{2 i m \pi },-k,\frac{-\frac{\pi  (2 j-\beta
   +1)}{\beta }-i \log (a)+\pi }{2 \pi }\right)\right)\\
=\frac{1}{2\pi i}\int_{C}\sum_{j=0}^{\beta-1}\frac{\pi  a^w w^{-k-1} \csc (\pi  (m+w)) \cos \left(\frac{\pi  (-\beta +2 j+1) (-\alpha
   +m+w)}{\beta }\right)}{2 \beta  \left(\cos \left(\frac{\pi  \alpha  (2 j+1)}{\beta }\right)+1\right)}dw
\end{multline}
Proceeding as we did before we next derive the forth contour integral representation in equation (\ref{eq2}) by using equation (\ref{eq3}) and setting $b= \pi ,t= 0,x= \frac{\pi  (-\alpha +2 j+1)}{\alpha }$ multiplying both sides by $\frac{\pi }{2 \alpha  \left(\cos \left(\frac{\pi  \beta  (2 j+1)}{\alpha }\right)+1\right)}$ and simplifying to get;
\begin{multline}\label{eq7}
-\sum_{j=0}^{\alpha-1}\frac{2^{k-1} (i \pi )^{k+1} e^{i \left(\pi  m-\frac{\pi  m (-\alpha +2 j+1)}{\alpha }\right)} }{\alpha  \Gamma(k+1) \left(\cos \left(\frac{\pi  \beta  (2 j+1)}{\alpha
   }\right)+1\right)}\left(\Phi
   \left(e^{2 i m \pi },-k,\frac{-\frac{\pi  (2 j-\alpha +1)}{\alpha }-i \log (a)+\pi }{2 \pi }\right)\right. \\ \left.
+e^{\frac{2 i \pi
    m (-\alpha +2 j+1)}{\alpha }} \Phi \left(e^{2 i m \pi },-k,\frac{\frac{\pi  (2 j-\alpha +1)}{\alpha }-i \log
   (a)+\pi }{2 \pi }\right)\right)\\
=\frac{1}{2\pi i}\int_{C}\sum_{j=0}^{\alpha-1}\frac{\pi  a^w w^{-k-1} \csc (\pi  (m+w)) \cos \left(\frac{\pi  (-\alpha +2 j+1) (m+w)}{\alpha
   }\right)}{2 \alpha  \left(\cos \left(\frac{\pi  \beta  (2 j+1)}{\alpha }\right)+1\right)}dw
\end{multline}
We finally tackle writing down the final contour integral representation on the right-hand side of equation (\ref{eq2}) by using equation (\ref{eq3}) and setting $b= \pi ,t= -\beta ,x= \frac{\pi  (-\alpha +2 j+1)}{\alpha }$, multiplying both sides by $\frac{\pi }{2 \alpha  \left(\cos \left(\frac{\pi  \beta  (2 j+1)}{\alpha }\right)+1\right)}$ and simplifying to get;
\begin{multline}\label{eq8}
-\sum_{j=0}^{\alpha-1}\frac{2^{k-1} (i \pi )^{k+1} e^{i \left(\pi  m-\frac{\pi  (-\alpha +2 j+1) (m-\beta )}{\alpha }\right)}}{\alpha  \Gamma(k+1) \left(\cos \left(\frac{\pi  \beta  (2 j+1)}{\alpha }\right)+1\right)} \left(e^{\frac{2 i \pi  (-\alpha +2 j+1) (m-\beta )}{\alpha }} \Phi \left(e^{2 i m \pi },-k,\frac{\frac{\pi  (2 j-\alpha +1)}{\alpha }-i \log (a)+\pi }{2 \pi }\right)\right. \\ \left.
+\Phi \left(e^{2 i m \pi },-k,\frac{-\frac{\pi  (2 j-\alpha+1)}{\alpha }-i \log (a)+\pi }{2 \pi }\right)\right)\\
=\frac{1}{2\pi i}\int_{C}\sum_{j=0}^{\alpha-1}\frac{\pi  a^w w^{-k-1} \csc (\pi  (m+w)) \cos \left(\frac{\pi  (-\alpha +2 j+1) (-\beta+m+w)}{\alpha }\right)}{2 \alpha  \left(\cos \left(\frac{\pi  \beta  (2 j+1)}{\alpha }\right)+1\right)}dw
\end{multline}
It is worthwhile to note that the interchanging of the contour integral and definite integral and finite series are valid based on Tonelli's theorem for sums and integrals, see page 177 in \cite{gelca} as the summand and integral are of bounded measure over the space $\mathbb{C} \times [0,\infty)$, $\mathbb{C} \times [0,\alpha)$ and $\mathbb{C} \times [0,\beta)$.
\section{Main results}
The main theorem in this work is derived by observing that the addition of  right-hand sides of equations (\ref{eq5}), (\ref{eq6}), (\ref{eq7}) and (\ref{eq8}) is equal to the right-hand side of equation (\ref{eq4}) so we can equate the left-hand sides and simplify relative to equation (\ref{eq2}) to yield the stated theorem.
\begin{theorem}
\begin{multline}\label{eq9}
\int_0^{\infty } \frac{x^{-1+m} \log ^k(a x)}{\left(1+x^{\alpha }\right) \left(1+x^{\beta }\right)} \,
   dx\\
=\sum _{j=0}^{\beta -1} \frac{1}{\beta  \left(1+\cos \left(\frac{(1+2 j) \pi  \alpha }{\beta }\right)\right)}(2 i)^{-1+k} \exp \left(-\frac{i \pi  (m (2+4 j-3 \beta )+\alpha  (-1-2 j+\beta ))}{\beta }\right) \\
\left(e^{\frac{i m \pi  (1+2 j-\beta )}{\beta }}+e^{\frac{i \pi  (m-\alpha ) (1+2 j-\beta)}{\beta }}\right) \pi ^{1+k} \left(e^{\frac{i \pi  (2 m-\alpha ) (1+2 j-\beta )}{\beta }} \Phi \left(e^{2 i m \pi },-k,\frac{\pi +2 j \pi -i \beta  \log (a)}{2 \pi  \beta }\right)\right. \\ \left.
+\Phi \left(e^{2 i m \pi },-k,-\frac{\pi +2 j \pi -2 \pi  \beta +i \beta  \log (a)}{2 \pi  \beta }\right)\right)\\
+\sum _{j=0}^{\alpha -1} \frac{1}{\alpha  \left(1+\cos \left(\frac{(1+2 j) \pi  \beta }{\alpha }\right)\right)}(2 i)^{-1+k} \exp \left(-\frac{i \pi  (m(2+4 j-3 \alpha )+(-1-2 j+\alpha ) \beta )}{\alpha }\right)\\ \left(e^{\frac{i m \pi  (1+2 j-\alpha )}{\alpha }}+e^{\frac{i \pi  (1+2 j-\alpha ) (m-\beta )}{\alpha }}\right) \pi ^{1+k} \left(e^{\frac{i \pi  (1+2 j-\alpha ) (2 m-\beta )}{\alpha }} \Phi \left(e^{2 i m \pi },-k,\frac{\pi +2 j \pi -i \alpha  \log (a)}{2 \pi  \alpha}\right)\right. \\ \left.
+\Phi \left(e^{2 i m \pi },-k,-\frac{\pi +2 j \pi -2 \pi  \alpha +i \alpha  \log (a)}{2 \pi  \alpha }\right)\right)
\end{multline}
where $Re(m)>0,\alpha \neq \beta,\alpha,\beta$ are even.
\end{theorem}
\begin{example}
In this example we simply used equation (\ref{eq9}) and set $a=1$ and simplify.
\begin{multline}
\int_0^{\infty } \frac{x^{-1+m} \log ^k(x)}{\left(1+x^{\alpha }\right) \left(1+x^{\beta }\right)} \,
   dx\\
=\sum _{j=0}^{-1+\alpha } \frac{1}{\alpha  \left(1+\cos \left(\frac{(1+2 j) \pi  \beta
   }{\alpha }\right)\right)}(2 i)^{-1+k} \exp \left(-\frac{i \pi  (m (2+4 j-3 \alpha )+(-1-2 j+\alpha )
   \beta )}{\alpha }\right)\\
 \left(e^{\frac{i m \pi  (1+2 j-\alpha )}{\alpha }}+e^{\frac{i \pi  (1+2 j-\alpha )
   (m-\beta )}{\alpha }}\right) \pi ^{1+k} \left(e^{\frac{i \pi  (1+2 j-\alpha ) (2 m-\beta )}{\alpha }} \Phi
   \left(e^{2 i m \pi },-k,\frac{\pi +2 j \pi }{2 \pi  \alpha }\right)\right. \\ \left.
+\Phi \left(e^{2 i m \pi },-k,-\frac{\pi +2
   j \pi -2 \pi  \alpha }{2 \pi  \alpha }\right)\right)\\
+\sum _{j=0}^{-1+\beta } \frac{1}{\beta  \left(1+\cos \left(\frac{(1+2 j)
   \pi  \alpha }{\beta }\right)\right)}(2 i)^{-1+k} \exp \left(-\frac{i \pi  (m (2+4 j-3 \beta
   )+\alpha  (-1-2 j+\beta ))}{\beta }\right)\\
 \left(e^{\frac{i m \pi  (1+2 j-\beta )}{\beta }}+e^{\frac{i \pi 
   (m-\alpha ) (1+2 j-\beta )}{\beta }}\right) \pi ^{1+k} \left(e^{\frac{i \pi  (2 m-\alpha ) (1+2 j-\beta
   )}{\beta }} \Phi \left(e^{2 i m \pi },-k,\frac{\pi +2 j \pi }{2 \pi  \beta }\right)\right. \\ \left.
+\Phi \left(e^{2 i m \pi
   },-k,-\frac{\pi +2 j \pi -2 \pi  \beta }{2 \pi  \beta }\right)\right)
\end{multline}
\end{example}
\begin{example}
In the next example we used equation (\ref{eq9}) and set $a=1,m=1/2$, simplify in terms of the Hurwitz zeta function see [DLMF,\href{https://dlmf.nist.gov/25.11}{25.11}] and take the first partial derivative with respect to $k$ and set $k=0$ and simplify using [DLMF,\href{https://dlmf.nist.gov/25.11.E18}{25.11.18}].
\begin{multline}
\int_0^{\infty } \frac{\log (\log (x))}{\sqrt{x} \left(1+x^{\alpha }\right) \left(1+x^{\beta }\right)} \,
   dx\\
=\sum _{j=0}^{\beta -1} \frac{1}{8 \beta  \left(1+\cos \left(\frac{(1+2 j) \pi  \alpha }{\beta }\right)\right)}\left(e^{\frac{i \pi  (1+2 j-\beta )}{2 \beta }}+e^{\frac{i \pi 
   \left(\frac{1}{2}-\alpha \right) (1+2 j-\beta )}{\beta }}\right) \pi \\
 \left(i \pi +\log (16)+2 \log (\pi )-4 \log \left(-1+\frac{1+2 j}{4 \beta }\right)\right. \\ \left.
+4 \log \left(\frac{1+2 j-2 \beta }{4 \beta }\right)-4 \text{log$\Gamma $}\left(-1+\frac{1+2 j}{4 \beta }\right)+4 \text{log$\Gamma $}\left(\frac{1+2 j-2 \beta }{4\beta }\right)\right. \\ \left.
+e^{\frac{i \pi  (-1+\alpha ) (1+2 j-\beta )}{\beta }} \left(i \pi +\log (16)+2 \log (\pi )+4 \log \left(-\frac{1+2 j}{\beta }\right)\right.\right. \\ \left.\left.
-4 \log \left(-\frac{1+2 j+2 \beta }{\beta }\right)
+4 \text{log$\Gamma $}\left(-\frac{1+2 j}{4 \beta }\right)-4 \text{log$\Gamma $}\left(-\frac{1+2 j+2 \beta }{4 \beta }\right)\right)\right)\\
+\sum
   _{j=0}^{\alpha -1} \frac{1}{8 \alpha  \left(1+\cos \left(\frac{(1+2 j) \pi  \beta }{\alpha
   }\right)\right)}\left(e^{\frac{i \pi  (1+2 j-\alpha )}{2 \alpha }}+e^{\frac{i \pi  (1+2 j-\alpha )
   \left(\frac{1}{2}-\beta \right)}{\alpha }}\right) \pi \\
 \left(i \pi +\log (16)+2 \log (\pi )-4 \log
   \left(-1+\frac{1+2 j}{4 \alpha }\right)+4 \log \left(\frac{1+2 j-2 \alpha }{4 \alpha }\right)\right. \\ \left.
-4\text{log$\Gamma $}\left(-1+\frac{1+2 j}{4 \alpha }\right)+4 \text{log$\Gamma $}\left(\frac{1+2 j-2 \alpha }{4 \alpha }\right)\right. \\ \left.
+e^{\frac{i \pi  (1+2 j-\alpha ) (-1+\beta )}{\alpha }} \left(i \pi +\log (16)+2 \log (\pi )+4 \log \left(-\frac{1+2 j}{\alpha }\right)\right.\right. \\ \left.\left.
-4 \log \left(-\frac{1+2 j+2 \alpha }{\alpha }\right)+4
   \text{log$\Gamma $}\left(-\frac{1+2 j}{4 \alpha }\right)-4 \text{log$\Gamma $}\left(-\frac{1+2 j+2 \alpha }{4 \alpha }\right)\right)\right)
\end{multline}
\end{example}
\begin{example}
In the next example we used equation (\ref{eq9}) and set $m=1/4$, and take the first partial derivative with respect to $k$ and set $k=0$ and simplify using [AIMS,\href{https://doi.org/10.3934/math.20221021}{6.7}].
\begin{multline}
\int_0^{\infty } \frac{\log (\log (a x))}{x^{3/4} \left(1+x^{\alpha }\right) \left(1+x^{\beta }\right)} \,
   dx\\
=\sum _{j=0}^{\beta -1} \frac{1}{\beta  \left(1+\cos \left(\frac{(1+2 j) \pi  \alpha }{\beta
   }\right)\right)}\left(\frac{1}{8}+\frac{i}{8}\right) e^{\frac{i \pi  (-2+j (-4+8 \alpha )-4
   \alpha  (-1+\beta )+3 \beta )}{4 \beta }} \left(e^{\frac{i \pi  (1+2 j-\beta )}{4 \beta }}+e^{\frac{i \pi 
   \left(\frac{1}{4}-\alpha \right) (1+2 j-\beta )}{\beta }}\right) \pi \\
 \left(\pi -2 i \log (2 \pi )+e^{\frac{i \pi  \left(\frac{1}{2}-\alpha \right) (1+2 j-\beta )}{\beta }} \left(\pi -2 i \log (2 \pi )+(2+2 i) \left(\log \left(\frac{\Gamma \left(\frac{\pi +2 j \pi -i \beta  \log (a)}{8 \pi  \beta }\right)}{2 \Gamma\left(\frac{1}{4} \left(2+\frac{\pi +2 j \pi -i \beta  \log (a)}{2 \pi  \beta }\right)\right)}\right)\right.\right.\right. \\ \left.\left.\left.
+i \log\left(\frac{\Gamma \left(\frac{1}{4} \left(1+\frac{\pi +2 j \pi -i \beta  \log (a)}{2 \pi  \beta
   }\right)\right)}{2 \Gamma \left(\frac{1}{4} \left(3+\frac{\pi +2 j \pi -i \beta  \log (a)}{2 \pi  \beta
   }\right)\right)}\right)\right)\right)+(2+2 i) \left(\log \left(\frac{\Gamma \left(-\frac{\pi +2 j \pi -2 \pi 
   \beta +i \beta  \log (a)}{8 \pi  \beta }\right)}{2 \Gamma \left(\frac{1}{4} \left(2-\frac{\pi +2 j \pi -2 \pi 
   \beta +i \beta  \log (a)}{2 \pi  \beta }\right)\right)}\right)\right.\right. \\ \left.\left.
+i \log \left(\frac{\Gamma \left(\frac{1}{4}
   \left(1-\frac{\pi +2 j \pi -2 \pi  \beta +i \beta  \log (a)}{2 \pi  \beta }\right)\right)}{2 \Gamma
   \left(\frac{1}{4} \left(3-\frac{\pi +2 j \pi -2 \pi  \beta +i \beta  \log (a)}{2 \pi  \beta
   }\right)\right)}\right)\right)\right)\\
+\sum _{j=0}^{\alpha -1} \frac{1}{\alpha  \left(1+\cos \left(\frac{(1+2 j) \pi  \beta }{\alpha
   }\right)\right)}\left(\frac{1}{8}+\frac{i}{8}\right) e^{\frac{i \pi  (-2+\alpha 
   (3-4 \beta )+4 \beta +j (-4+8 \beta ))}{4 \alpha }} \left(e^{\frac{i \pi  (1+2 j-\alpha )}{4 \alpha
   }}+e^{\frac{i \pi  (1+2 j-\alpha ) \left(\frac{1}{4}-\beta \right)}{\alpha }}\right) \pi  \\
\left(\pi -2 i \log (2 \pi )+e^{\frac{i \pi  (1+2 j-\alpha ) \left(\frac{1}{2}-\beta \right)}{\alpha }} \left(\pi -2 i \log (2 \pi)+(2+2 i) \left(\log \left(\frac{\Gamma \left(\frac{\pi +2 j \pi -i \alpha  \log (a)}{8 \pi  \alpha }\right)}{2
   \Gamma \left(\frac{1}{4} \left(2+\frac{\pi +2 j \pi -i \alpha  \log (a)}{2 \pi  \alpha }\right)\right)}\right)\right.\right.\right. \\ \left.\left.\left.
+i \log \left(\frac{\Gamma \left(\frac{1}{4} \left(1+\frac{\pi +2 j \pi -i \alpha  \log (a)}{2 \pi  \alpha
   }\right)\right)}{2 \Gamma \left(\frac{1}{4} \left(3+\frac{\pi +2 j \pi -i \alpha  \log (a)}{2 \pi  \alpha
   }\right)\right)}\right)\right)\right)+(2+2 i) \left(\log \left(\frac{\Gamma \left(-\frac{\pi +2 j \pi -2 \pi 
   \alpha +i \alpha  \log (a)}{8 \pi  \alpha }\right)}{2 \Gamma \left(\frac{1}{4} \left(2-\frac{\pi +2 j \pi -2 \pi
    \alpha +i \alpha  \log (a)}{2 \pi  \alpha }\right)\right)}\right)\right.\right. \\ \left.\left.
+i \log \left(\frac{\Gamma \left(\frac{1}{4}
   \left(1-\frac{\pi +2 j \pi -2 \pi  \alpha +i \alpha  \log (a)}{2 \pi  \alpha }\right)\right)}{2 \Gamma
   \left(\frac{1}{4} \left(3-\frac{\pi +2 j \pi -2 \pi  \alpha +i \alpha  \log (a)}{2 \pi  \alpha
   }\right)\right)}\right)\right)\right)
\end{multline}
\end{example}
\begin{example}
In the next example we used equation (\ref{eq9}) and set $m=1$ and simplify using equation [DLMF,\href{https://dlmf.nist.gov/25.14.E2}{25.14.2}],
\begin{multline}
\int_0^{\infty } \frac{\log ^k(a x)}{\left(1+x^{\alpha }\right) \left(1+x^{\beta }\right)} \, dx\\
=\sum_{j=0}^{\beta -1} \frac{1}{\beta  \left(1+\cos \left(\frac{(1+2 j) \pi  \alpha }{\beta }\right)\right)}(2 i)^{-1+k} \exp \left(\frac{i \pi  (1-2 j (-1+\alpha )+\alpha  (-1+\beta ))}{\beta
   }\right) \\
\left(1+e^{\frac{i \pi  \alpha  (1+2 j-\beta )}{\beta }}\right) \pi ^{1+k} \left(\zeta
   \left(-k,\frac{\pi +2 j \pi -i \beta  \log (a)}{2 \pi  \beta }\right)\right. \\ \left.+e^{\frac{i \pi  (-2+\alpha ) (1+2 j-\beta
   )}{\beta }} \zeta \left(-k,-\frac{\pi +2 j \pi -2 \pi  \beta +i \beta  \log (a)}{2 \pi  \beta
   }\right)\right)\\
+\sum _{j=0}^{\alpha -1} \frac{1}{\alpha  \left(1+\cos \left(\frac{(1+2 j) \pi  \beta }{\alpha }\right)\right)}(2 i)^{-1+k} \exp \left(\frac{i \pi  (1-2 j (-1+\beta )+(-1+\alpha ) \beta )}{\alpha}\right)\\
 \left(1+e^{\frac{i \pi  (1+2 j-\alpha ) \beta }{\alpha }}\right) \pi ^{1+k} \left(\zeta \left(-k,\frac{\pi +2 j \pi -i \alpha  \log (a)}{2 \pi  \alpha }\right)\right. \\ \left.+e^{\frac{i \pi  (1+2 j-\alpha )(-2+\beta )}{\alpha }} \zeta \left(-k,-\frac{\pi +2 j \pi -2 \pi  \alpha +i \alpha  \log (a)}{2 \pi  \alpha}\right)\right)
\end{multline}
\end{example}
\begin{example}
In the next example we used equation (\ref{eq9}) and set $m=1$ and simplify using equation [DLMF,\href{https://dlmf.nist.gov/25.14.E2}{25.14.2}], next we take the first partial derivative with respect to $k$ and set $k=0$ and simplify using [DLMF,\href{https://dlmf.nist.gov/25.11.E18}{25.11.18}].
\begin{multline}\label{eq26}
\int_0^{\infty } \frac{\log (\log (a x))}{\left(1+x^{\alpha }\right) \left(1+x^{\beta }\right)} \, dx\\
=\sum_{j=0}^{\beta -1} \frac{1}{8 \beta ^2 \left(1+\cos \left(\frac{(1+2 j) \pi  \alpha }{\beta
   }\right)\right)}\exp \left(\frac{i \pi  (1-2 j (-1+\alpha )+\alpha  (-1+\beta ))}{\beta }\right)
   \left(1+e^{\frac{i \pi  \alpha  (1+2 j-\beta )}{\beta }}\right)\\
 \left(e^{\frac{i \pi  (-2+\alpha ) (1+2 j-\beta
   )}{\beta }} (\pi +2 j \pi -\pi  \beta +i \beta  \log (a)) (\pi -2 i (\log (2)+\log (\pi )))+(\pi  (-1-2 j+\beta
   )\right. \\ \left.
+i \beta  \log (a)) (\pi -2 i (\log (2)+\log (\pi )))+4 i \pi  \beta  \left(-\frac{1}{2} \log (2 \pi )+\log
   \left(-1+\frac{\pi +2 j \pi -i \beta  \log (a)}{2 \pi  \beta }\right)\right.\right. \\ \left.\left.
+\text{log$\Gamma $}\left(-1+\frac{\pi +2 j \pi -i \beta  \log (a)}{2 \pi  \beta }\right)\right)\right. \\ \left.
+4 i e^{\frac{i \pi  (-2+\alpha ) (1+2 j-\beta )}{\beta }}\pi  \beta  \left(-\frac{1}{2} \log (2 \pi )+\log \left(-1-\frac{\pi +2 j \pi -2 \pi  \beta +i \beta  \log (a)}{2 \pi  \beta }\right)\right.\right. \\ \left.\left.
+\text{log$\Gamma $}\left(-1-\frac{\pi +2 j \pi -2 \pi  \beta +i \beta  \log (a)}{2\pi  \beta }\right)\right)\right)\\
+\sum _{j=0}^{\alpha -1} \frac{1}{8 \alpha^2 \left(1+\cos \left(\frac{(1+2 j) \pi  \beta }{\alpha }\right)\right)}
\exp \left(\frac{i \pi  (1-2 j (-1+\beta )+(-1+\alpha ) \beta )}{\alpha }\right) \left(1+e^{\frac{i \pi  (1+2 j-\alpha ) \beta }{\alpha }}\right)\\
 \left(e^{\frac{i \pi  (1+2 j-\alpha ) (-2+\beta )}{\alpha }} (\pi +2 j \pi -\pi  \alpha +i \alpha  \log (a)) (\pi -2 i (\log (2)+\log (\pi)))+(\pi  (-1-2 j+\alpha )\right. \\ \left.
+i \alpha  \log (a)) (\pi -2 i (\log (2)+\log (\pi )))+4 i \pi  \alpha  \left(-\frac{1}{2} \log (2 \pi )+\log \left(-1+\frac{\pi +2 j \pi -i \alpha  \log (a)}{2 \pi  \alpha}\right)\right.\right. \\ \left.\left.
+\text{log$\Gamma $}\left(-1+\frac{\pi +2 j \pi -i \alpha  \log (a)}{2 \pi  \alpha }\right)\right)\right. \\ \left.+4 i e^{\frac{i \pi  (1+2 j-\alpha ) (-2+\beta )}{\alpha }} \pi  \alpha  \left(-\frac{1}{2} \log (2 \pi )+\log \left(-1-\frac{\pi +2 j \pi -2 \pi  \alpha +i \alpha  \log (a)}{2 \pi  \alpha }\right)\right.\right. \\ \left.\left.
+\text{log$\Gamma
   $}\left(-1-\frac{\pi +2 j \pi -2 \pi  \alpha +i \alpha  \log (a)}{2 \pi  \alpha }\right)\right)\right)
\end{multline}
\end{example}
\begin{example}
In this next example we used equation (\ref{eq26}) and set $a=1,\alpha=2,\beta=4$ and simplify,
\begin{multline}
\int_0^{\infty } \frac{\log (\log (x))}{\left(1+x^2\right) \left(1+x^4\right)} \, dx=\frac{1}{2} \pi  \log
   \left(\frac{\left(\frac{1}{3}+\frac{i}{3}\right) \left(1+\sqrt{2}\right)^{\frac{i}{\sqrt{2}}} \sqrt{\pi } \Gamma
   \left(-\frac{1}{4}\right)}{\Gamma \left(-\frac{3}{4}\right)}\right)
\end{multline}
\end{example}
\begin{example}
In this next example we used equation (\ref{eq26}) and set $a=1,\alpha=2,\beta=8$ and simplify,
\begin{multline}
\int_0^{\infty } \frac{\log (\log (x))}{\left(1+x^2\right) \left(1+x^8\right)} \, dx\\
=-\frac{1}{8}
   (-1)^{5/8} \pi  \log \left(4^{(-1)^{3/8}} e^{(-1)^{7/8} \pi } \pi ^{2 (-1)^{3/8}} \cot
   ^{(-1+i)-\sqrt{2}}\left(\frac{3 \pi }{16}\right) \left(\frac{3 \Gamma \left(-\frac{3}{4}\right)}{\Gamma
   \left(-\frac{1}{4}\right)}\right)^{-4 (-1)^{3/8}}\right. \\ \left. \tan ^{(1+i)-i \sqrt{2}}\left(\frac{\pi}{16}\right)\right)
\end{multline}
\end{example}
\begin{example}
This next example involves a bit of algebraic manipulation of some very long equations. We started by using equation (\ref{eq9}) and replacing $m$ by $u+iv$. Next we replaced $v$ by $-v$ and took the difference of the two last equations formulated.
\begin{multline}\label{eq29}
\int_0^{\infty } \frac{x^{-1+u-i v} \left(1-x^{2 i v}\right)}{\left(1+x^{\alpha }\right) \left(1+x^{\beta
   }\right) \log (x)} \, dx\\
=\sum _{j=0}^{\beta -1} \frac{i}{4 \beta  \left(1+\cos \left(\frac{(1+2 j)
   \pi  \alpha }{\beta }\right)\right)}\ \exp \left(-\frac{1}{2} \pi  \left(i+4 v+2 i \alpha +\frac{2
   i (1+2 j) (u-i v+\alpha )}{\beta }\right)\right) \\
\left(e^{i \pi  \alpha }+e^{\frac{i (1+2 j) \pi  \alpha }{\beta
   }}\right) \left(\exp \left(\frac{i \pi  ((2+4 j) u+(-2 i v+\alpha ) \beta )}{\beta }\right) \Phi \left(e^{2 i \pi 
   (u+i v)},1,\frac{\pi +2 j \pi }{2 \pi  \beta }\right)\right. \\ \left.
+e^{2 i \pi  u+\frac{(1+2 j) \pi  (2 v+i \alpha )}{\beta }}
   \Phi \left(e^{2 i \pi  (u+i v)},1,-\frac{\pi +2 j \pi -2 \pi  \beta }{2 \pi  \beta }\right)\right. \\ \left.
-\exp \left(\pi  \left(2v+i \alpha +\frac{2 (1+2 j) (i u+v)}{\beta }\right)\right) \Phi \left(e^{2 \pi  (i u+v)},1,\frac{\pi +2 j \pi }{2 \pi  \beta }\right)\right. \\ \left.
-e^{\frac{i \pi  (\alpha +2 j \alpha +2 u \beta -4 i v \beta )}{\beta }} \Phi \left(e^{2 \pi  (iu+v)},1,-\frac{\pi +2 j \pi -2 \pi  \beta }{2 \pi  \beta }\right)\right)\\
+\sum _{j=0}^{\alpha -1} \frac{i }{2 \left(1+e^{\frac{i (1+2 j) \pi  \beta }{\alpha }}\right)^2 \alpha }\exp \left(-\frac{\pi  (2 i (1+2 j) u+v (2+4 j+4 \alpha )-i \alpha  (-1-2 \beta ))}{2 \alpha }\right)\\
 \left(e^{i \pi  \beta }+e^{\frac{i (1+2 j) \pi  \beta }{\alpha }}\right)
 \left(\exp \left(\frac{i \pi  ((2+4 j) u+\alpha  (-2 i v+\beta ))}{\alpha }\right) \Phi\left(e^{2 i \pi  (u+i v)},1,\frac{\pi +2 j \pi }{2 \pi  \alpha }\right)\right. \\ \left.
+\exp \left(\frac{\pi  ((2+4 j) v+i (2 u\alpha +\beta +2 j \beta ))}{\alpha }\right) \Phi \left(e^{2 i \pi  (u+i v)},1,-\frac{\pi +2 j \pi -2 \pi  \alpha }{2 \pi  \alpha }\right)\right. \\ \left.
-\exp \left(\frac{\pi  (2 i (1+2 j) u+2 v (1+2 j+\alpha )+i \alpha  \beta )}{\alpha }\right) \Phi \left(e^{2 \pi  (i u+v)},1,\frac{\pi +2 j \pi }{2 \pi  \alpha }\right)\right. \\ \left.
-e^{\frac{i \pi  (2 u \alpha -4 i v \alpha +\beta +2 j \beta )}{\alpha }} \Phi \left(e^{2 \pi  (i u+v)},1,-\frac{\pi +2 j \pi -2 \pi  \alpha }{2 \pi  \alpha }\right)\right)
\end{multline}
where $Re(u)\neq 0.$
\end{example}
\begin{example}
In this example we used equation (\ref{eq9}) and set $m=1/2$ then simplified in terms of the Hurwitz zeta function using entry (4) in the table below (64:12:5) in \cite{atlas}. Next we take the limit of both sides and apply l'Hopital's rule to the right-hand side as $k\to -1$ simplify using [Wolfram,\href{https://mathworld.wolfram.com/HurwitzZetaFunction.html}{19}] and replace $a$ by $e^{a\pi}$.
\begin{multline}\label{eq210}
\int_0^{\infty } \frac{1}{\sqrt{x} \left(1+x^{\alpha }\right) \left(1+x^{\beta }\right) (a \pi +\log (x))} \,
   dx\\
=\sum _{j=0}^{\beta -1} \frac{i }{8 \beta  \left(1+\cos \left(\frac{(1+2 j) \pi  \alpha }{\beta }\right)\right)}\left(e^{\frac{i \pi  (1+2 j-\beta )}{2 \beta }}+e^{\frac{i \pi 
   \left(\frac{1}{2}-\alpha \right) (1+2 j-\beta )}{\beta }}\right) \\\left(\psi ^{(0)}\left(\frac{1+2 j-i a \beta }{4
   \beta }\right)-\psi ^{(0)}\left(\frac{1+2 j+2 \beta -i a \beta }{4 \beta }\right)\right. \\ \left.
+e^{\frac{i \pi  (-1+\alpha ) (1+2j-\beta )}{\beta }} \left(-\psi ^{(0)}\left(-\frac{1+2 j-4 \beta +i a \beta }{4 \beta }\right)+\psi
   ^{(0)}\left(-\frac{1+2 j-2 \beta +i a \beta }{4 \beta }\right)\right)\right)\\
+\sum _{j=0}^{\alpha -1} \frac{i }{8\alpha  \left(1+\cos \left(\frac{(1+2 j) \pi  \beta }{\alpha }\right)\right)}\left(e^{\frac{i \pi  (1+2 j-\alpha )}{2 \alpha }}+e^{\frac{i \pi  (1+2 j-\alpha ) \left(\frac{1}{2}-\beta \right)}{\alpha }}\right) \\\left(\psi^{(0)}\left(\frac{1+2 j-i a \alpha }{4 \alpha }\right)-\psi ^{(0)}\left(\frac{1+2 j+2 \alpha -i a \alpha }{4 \alpha }\right)\right. \\ \left.
+e^{\frac{i \pi  (1+2 j-\alpha ) (-1+\beta )}{\alpha }} \left(-\psi ^{(0)}\left(-\frac{1+2 j-4 \alpha +i a \alpha }{4 \alpha }\right)+\psi ^{(0)}\left(-\frac{1+2 j-2 \alpha +i a \alpha }{4 \alpha }\right)\right)\right)
\end{multline}
where $Re(a) \geq 0.$
\end{example}
\begin{example}
In the following example we used equation (\ref{eq210}) simply replaced $a$ by $-a$ and took the difference of the two equations.
\begin{multline}
\int_0^{\infty } \frac{2 a \pi }{\sqrt{x} \left(1+x^{\alpha }\right) \left(1+x^{\beta }\right) \left(a^2 \pi
   ^2-\log ^2(x)\right)} \, dx\\
=\sum _{j=0}^{\beta -1} \frac{i }{8 \beta  \left(1+\cos \left(\frac{(1+2 j) \pi  \alpha }{\beta
   }\right)\right)}\left(e^{\frac{i \pi  (1+2 j-\beta )}{2 \beta}}+e^{\frac{i \pi  \left(\frac{1}{2}-\alpha \right) (1+2 j-\beta )}{\beta }}\right)\\
 \left(\psi ^{(0)}\left(\frac{1+2 j-i a \beta }{4 \beta }\right)+e^{\frac{i \pi  (-1+\alpha ) (1+2 j-\beta )}{\beta }} \psi ^{(0)}\left(-\frac{1+2 j-4\beta -i a \beta }{4 \beta }\right)\right. \\ \left.
-e^{\frac{i \pi  (-1+\alpha ) (1+2 j-\beta )}{\beta }} \psi^{(0)}\left(-\frac{1+2 j-2 \beta -i a \beta }{4 \beta }\right)-\psi ^{(0)}\left(\frac{1+2 j+2 \beta -i a \beta }{4 \beta }\right)\right. \\ \left.
-\psi ^{(0)}\left(\frac{1+2 j+i a \beta }{4 \beta }\right)-e^{\frac{i \pi  (-1+\alpha ) (1+2 j-\beta )}{\beta }} \psi ^{(0)}\left(-\frac{1+2 j-4 \beta +i a \beta }{4 \beta }\right)\right. \\ \left.
+e^{\frac{i \pi  (-1+\alpha ) (1+2 j-\beta )}{\beta }} \psi ^{(0)}\left(-\frac{1+2 j-2 \beta +i a \beta }{4 \beta }\right)+\psi ^{(0)}\left(\frac{1+2 j+2 \beta +i a \beta }{4 \beta }\right)\right)\\
+\sum _{j=0}^{\alpha -1} \frac{i }{8 \alpha  \left(1+\cos \left(\frac{(1+2 j) \pi  \beta }{\alpha
   }\right)\right)}\left(e^{\frac{i \pi  (1+2 j-\alpha )}{2 \alpha }}+e^{\frac{i \pi 
   (1+2 j-\alpha ) \left(\frac{1}{2}-\beta \right)}{\alpha }}\right)\\
 \left(\psi ^{(0)}\left(\frac{1+2 j-i a \alpha }{4\alpha }\right)+e^{\frac{i \pi  (1+2 j-\alpha ) (-1+\beta )}{\alpha }} \psi ^{(0)}\left(-\frac{1+2 j-4 \alpha -i a\alpha }{4 \alpha }\right)\right. \\ \left.
-e^{\frac{i \pi  (1+2 j-\alpha ) (-1+\beta )}{\alpha }} \psi ^{(0)}\left(-\frac{1+2 j-2 \alpha -i a \alpha }{4 \alpha }\right)-\psi ^{(0)}\left(\frac{1+2 j+2 \alpha -i a \alpha }{4 \alpha }\right)\right. \\ \left.
-\psi ^{(0)}\left(\frac{1+2 j+i a \alpha }{4 \alpha }\right)-e^{\frac{i \pi  (1+2 j-\alpha ) (-1+\beta )}{\alpha }} \psi^{(0)}\left(-\frac{1+2 j-4 \alpha +i a \alpha }{4 \alpha }\right)\right. \\ \left.
+e^{\frac{i \pi  (1+2 j-\alpha ) (-1+\beta )}{\alpha }} \psi ^{(0)}\left(-\frac{1+2 j-2 \alpha +i a \alpha }{4 \alpha }\right)+\psi ^{(0)}\left(\frac{1+2 j+2\alpha +i a \alpha }{4 \alpha }\right)\right)
\end{multline}
where $Re(a)>0,Im(a)> 0$.
\end{example}
\begin{example}
In this example we used equation (\ref{eq26}) and set $a=1,\alpha=4,\beta=6$ and simplify,
\begin{multline}
\int_0^{\infty } \frac{\log (\log (x))}{\left(1+x^4\right) \left(1+x^6\right)} \, dx\\
=\log \left(3^{\frac{1}{12}
   \left(-2+3 \sqrt{2}\right) \pi } 5^{-\frac{1}{12} \left(2+\sqrt{6+12 i \sqrt{6}}\right) \pi } 7^{-\frac{1}{12}
   \left(-2+\sqrt{6+12 i \sqrt{6}}\right) \pi } e^{\frac{1}{24} i \left(-1+3 \sqrt{2}\right) \pi ^2} (2 \pi
   )^{\frac{1}{12} \left(-1+3 \sqrt{2}\right) \pi }\right. \\ \left.
 \left(\frac{11 \Gamma \left(-\frac{11}{12}\right)}{\Gamma
   \left(-\frac{5}{12}\right)}\right)^{\frac{1}{6} \left(\pi +i \sqrt{3} \pi \right)} \left(\frac{\Gamma
   \left(-\frac{1}{4}\right)}{\Gamma \left(-\frac{3}{4}\right)}\right)^{\pi /6} \left(\frac{\Gamma
   \left(-\frac{3}{8}\right) \Gamma \left(-\frac{1}{8}\right)}{\Gamma \left(-\frac{7}{8}\right) \Gamma
   \left(-\frac{5}{8}\right)}\right)^{\frac{\pi }{2 \sqrt{2}}} \left(\frac{\Gamma \left(-\frac{7}{12}\right)}{\Gamma
   \left(-\frac{1}{12}\right)}\right)^{\frac{1}{6} \left(\pi -i \sqrt{3} \pi \right)}\right)
\end{multline}
\end{example}
\begin{example}
In this example we use equation (\ref{eq9}) and set $a=1,m-3/2$ and simplify in terms of the Hurwitz eta function using [DLMF,\href{https://dlmf.nist.gov/25.14.E2}{25.14.2}]. Next we take the limit and apply l'Hopital's rule to the right-hand side as $k\to -1$ and simplify using [Wolfram,\href{https://mathworld.wolfram.com/HurwitzZetaFunction.html}{20}].
\begin{multline}\label{eq213}
\int_0^{\infty } \frac{x-1}{\sqrt{x} \left(1+x^{\alpha }\right) \left(1+x^{\beta }\right) \log (x)} \, dx\\
=\sum
   _{j=0}^{\beta -1} \frac{i}{16 \beta } \left(e^{\frac{i \pi  (1+j (2-4 \alpha )-2 \alpha +3 \beta )}{2 \beta }}
   \left(-1+e^{\frac{i (\pi +2 j \pi )}{\beta }}\right) \left(e^{i \pi  \alpha }+e^{\frac{i (1+2 j) \pi  \alpha }{\beta
   }}\right) \psi ^{(0)}\left(\frac{1+2 j}{4 \beta }\right)\right. \\ \left.
+e^{\frac{i \pi  (-1-2 j+\beta )}{2 \beta }}
   \left(1+e^{\frac{i \pi  (-1+\alpha ) (1+2 j-\beta )}{\beta }}+e^{\frac{i \pi  \alpha  (1+2 j-\beta )}{\beta
   }}+e^{\frac{i \pi  (-1-2 j+\beta )}{\beta }}\right) \psi ^{(0)}\left(-\frac{1+2 j-4 \beta }{4 \beta
   }\right)\right. \\ \left.
+e^{-\frac{i \pi  (3+6 j+(-3+2 \alpha ) \beta )}{2 \beta }} \left(-1+e^{\frac{i (\pi +2 j \pi )}{\beta
   }}\right) \left(e^{i \pi  \alpha }+e^{\frac{i (1+2 j) \pi  \alpha }{\beta }}\right) \psi ^{(0)}\left(-\frac{1+2 j-2
   \beta }{4 \beta }\right)\right. \\ \left.
-e^{\frac{i \pi  (1+j (2-4 \alpha )-2 \alpha +3 \beta )}{2 \beta }} \left(-1+e^{\frac{i (\pi
   +2 j \pi )}{\beta }}\right) \left(e^{i \pi  \alpha }+e^{\frac{i (1+2 j) \pi  \alpha }{\beta }}\right) \psi
   ^{(0)}\left(\frac{1+2 j+2 \beta }{4 \beta }\right)\right) \sec ^2\left(\frac{(1+2 j) \pi  \alpha }{2 \beta
   }\right)\\
+\sum _{j=0}^{\alpha -1} \frac{i }{16 \alpha }\left(e^{\frac{i \pi  (1+3 \alpha +j (2-4 \beta )-2 \beta )}{2
   \alpha }} \left(-1+e^{\frac{i (\pi +2 j \pi )}{\alpha }}\right) \left(e^{i \pi  \beta }+e^{\frac{i (1+2 j) \pi 
   \beta }{\alpha }}\right) \psi ^{(0)}\left(\frac{1+2 j}{4 \alpha }\right)\right. \\ \left.
+e^{\frac{i \pi  (-1-2 j+\alpha )}{2 \alpha
   }} \left(1+e^{\frac{i \pi  (-1-2 j+\alpha )}{\alpha }}+e^{\frac{i \pi  (1+2 j-\alpha ) (-1+\beta )}{\alpha
   }}+e^{\frac{i \pi  (1+2 j-\alpha ) \beta }{\alpha }}\right) \psi ^{(0)}\left(-\frac{1+2 j-4 \alpha }{4 \alpha
   }\right)\right. \\ \left.
+e^{-\frac{i \pi  (3+6 j+\alpha  (-3+2 \beta ))}{2 \alpha }} \left(-1+e^{\frac{i (\pi +2 j \pi )}{\alpha
   }}\right) \left(e^{i \pi  \beta }+e^{\frac{i (1+2 j) \pi  \beta }{\alpha }}\right) \psi ^{(0)}\left(-\frac{1+2 j-2
   \alpha }{4 \alpha }\right)\right. \\ \left.
-e^{\frac{i \pi  (1+3 \alpha +j (2-4 \beta )-2 \beta )}{2 \alpha }} \left(-1+e^{\frac{i
   (\pi +2 j \pi )}{\alpha }}\right) \left(e^{i \pi  \beta }+e^{\frac{i (1+2 j) \pi  \beta }{\alpha }}\right) \psi
   ^{(0)}\left(\frac{1+2 j+2 \alpha }{4 \alpha }\right)\right) \sec ^2\left(\frac{(1+2 j) \pi  \beta }{2 \alpha
   }\right)
\end{multline}
\end{example}
\begin{example}
In this example we used equation (\ref{eq213}) and set $\alpha=2,\beta=4$ and simplified.
\begin{equation}
\int_0^{\infty } \frac{x-1}{\sqrt{x} \left(1+x^2\right) \left(1+x^4\right) \log (x)} \, dx=\log \left(\cot
   \left(\frac{\pi }{8}\right)\right)
\end{equation}
\end{example}
\begin{example}
In this example we use equation (\ref{eq9}) and form two equations by setting $a=1,m=1/2$ for the first equation and $a=1,m=3/2$ for the second equation and simplify in terms of the Hurwitz zeta function. Next we take the difference of these two equations. Next we take the first partial derivative with respect to $k$ and set $k=0$ and simplify using [DLMF,\href{https://dlmf.nist.gov/25.11.E18}{25.11.18}].
\begin{multline}
\int_0^{\infty } \frac{(x-1) \log (\log (x))}{\sqrt{x} \left(1+x^2\right) \left(1+x^4\right)} \,
   dx=\left(\frac{1}{8}+\frac{i}{8}\right) \sqrt[8]{-1} \pi  \left(-i \left((-1-i)+\sqrt{2}\right) \pi -2
   \sqrt[8]{-1} \left(2+\sqrt[8]{-1}\right)\right. \\ \left.
 \log \left(\frac{7}{5}\right)+(4+4 i) \log (2)-(6-2 i) \log (3)-2
   (-1)^{3/4} \log \left(\frac{143}{15}\right)\right. \\ \left.
-\sqrt{2} \log (16)+\sqrt[8]{-1} \log (81)+\log \left(49\ 3^{6
   \sqrt[4]{-1}} 5^{-2+2 i} 143^{-2 i}\right)\right. \\ \left.
-2 \left(\left((-1-i)+\sqrt{2}\right) \log (\pi
   )-\left(-1+\sqrt[4]{-1}\right) \log \left(\Gamma \left(-\frac{15}{16}\right)\right)+2 \sqrt[8]{-1} \log
   \left(\Gamma \left(-\frac{7}{8}\right)\right)\right.\right. \\ \left.\left.
+\left(i+(-1)^{3/4}\right) \log \left(\Gamma
   \left(-\frac{13}{16}\right)\right)+\left(i+(-1)^{3/4}\right) \log \left(\Gamma
   \left(-\frac{11}{16}\right)\right)\right.\right. \\ \left.\left.
-2 \sqrt[8]{-1} \log \left(\Gamma
   \left(-\frac{5}{8}\right)\right)-\left(-1+\sqrt[4]{-1}\right) \log \left(\Gamma
   \left(-\frac{9}{16}\right)\right)+\left(-1+\sqrt[4]{-1}\right) \log \left(\Gamma
   \left(-\frac{7}{16}\right)\right)\right.\right. \\ \left.\left.
-2 \sqrt[8]{-1} \log \left(\Gamma
   \left(-\frac{3}{8}\right)\right)-\left(i+(-1)^{3/4}\right) \log \left(\Gamma
   \left(-\frac{5}{16}\right)\right)-\left(i+(-1)^{3/4}\right) \log \left(\Gamma
   \left(-\frac{3}{16}\right)\right)\right.\right. \\ \left.\left.
+2 \sqrt[8]{-1} \log \left(\Gamma
   \left(-\frac{1}{8}\right)\right)+\left(-1+\sqrt[4]{-1}\right) \log \left(\Gamma
   \left(-\frac{1}{16}\right)\right)\right)\right)
\end{multline}
\end{example}
%
%
%
%
%
%
%
\section{Conclusion}
In this paper, we have presented a method for deriving an integral transform along with some definite integrals using contour integration. Other definite integrals will be pursued using this contour integral method. The results presented were numerically verified for both real and imaginary and complex values of the parameters in the integrals using Mathematica by Wolfram.
\end{document}